\documentstyle[amscd,amssymb,epsf,11pt]{amsart}
\newtheorem{thm}{Theorem}[section]

\newtheorem{conj}[thm]{Conjecture}
\newtheorem{prop}[thm]{Proposition}

\theoremstyle{definition}
\newtheorem{defn}[thm]{Definition}

\newtheorem{rem}[thm]{Remark}

\def\H{{\frak H}}
\def\C{{\Bbb C}}

\def\P{{\Bbb P}}

\def\Z{{\Bbb Z}}
\def\Q{{\Bbb Q}}

\makeatletter
\def\theequation{\@arabic\c@equation}
\makeatother
\numberwithin{equation}{thm}

\def\spec{{\hbox{\rom{Spec}\,}}}

\def\mmm{{m}}

\def\Q{{\Bbb Q}}
\def\Z{{\Bbb Z}}

\def\P{{\Bbb P}}
\def\C{{\Bbb C}}
\def\I{{\cal I}}

\def\oo{\OO}

\def\whsq{\vbox to 5.8pt
{\offinterlineskip\hrule
\hbox to 5.8pt{\vrule height
5.1pt\hss\vrule height 5.1pt}\hrule}}

\def\oo{{\cal O}}

\def\H{{\cal H}}
\def\C{{\Bbb C}}

\def\P{{\Bbb P}}

\def\Z{{\Bbb Z}}
\def\Q{{\Bbb Q}}

\def\phi{\varphi}

\def\({\left(}
\def\){\right)}

\def\part{P(n)}

\def\<{\langle}
\def\>{\rangle}

\begin{document}

\title[A conjectural generating function for numbers of curves]
{A conjectural generating function for numbers of curves on 
surfaces}
 \keywords{Severi degrees, Gromov-Witten invariants, nodal curves,
modular forms}
\author{Lothar G\"ottsche}
\address{International Center for Theoretical Physics, Strada
Costiera 11, P.O. Box 586, 34100 Trieste, Italy }
\email{gottsche@@ictp.trieste.it}

\begin{abstract}
I give a conjectural generating function for the numbers of $\delta$-nodal
curves   in a linear system
of dimension $\delta$ on an algebraic surface. It reproduces  the results
of Vainsencher \cite{V2} for the case $\delta\le 6$ and Kleiman-Piene 
\cite{K-P} for the case $\delta\le 8$.  The numbers of curves are expressed 
in terms of five  universal power series, three of which I give explicitly
as quasimodular forms. This gives in particular the numbers of curves of
 arbitrary genus on a  K3 surface and an abelian
surface in terms of quasimodular forms, generalizing the formula of 
Yau-Zaslow for rational curves on K3 surfaces.
The coefficients of the other two power series
can be determined by comparing with the recursive formulas of Caporaso-Harris
for the  Severi degrees in $\P_2$. We verify the conjecture for  
genus $2$ curves on an abelian surface. We also discuss a link of this 
problem with Hilbert schemes of points.
\end{abstract}

 \maketitle

\section{Introduction}

Let $L$ be a line bundle on a projective algebraic surface $S$.
In the case $\delta\le 6$ Vainsencher \cite{V2} proved  formulas
for the numbers $t^S_\delta(L)$ of $\delta$-nodal curves in a general
 $\delta$-dimensional sub-linear system   of $|L|$. 
By a refining Vainsechers approach Kleiman-Piene
\cite{K-P} extended the results to $\delta\le 8$.
The formulas hold under the assumption
that  $L$ is a sufficiently high power of a very ample line bundle.

In this paper we want to give a conjectural generating function for the
numbers $t^S_\delta(L)$. We will have only partial success:
We are able to express the $t^S_\delta(L)$ in terms of five universal 
generating functions in one variable $q$.
Three of these are Fourier developments of 
  (quasi-)modular forms, the other two we
have not been able to identify: the formulas of \cite{C-H}
for the Severi degrees on $\P_2$ yield an algorithm for computing 
their coefficients and I computed them up to degree 28.   As the  functions 
 are universal, one would hope that there    exists a nice closed
expression for them. 

If the canonical divisor of the surface $S$ is numerically trivial,
only the quasimodular forms  appear in the generating function.
Thus we obtain (conjecturally) the numbers of curves of arbitrary genus on
a K3 surface and on an abelian surface as the Fourier coefficients
of quasimodular forms. The formulas generalize the  calculation
of \cite{Y-Z} for the numbers of rational curves on K3 surfaces.
The fact that for K3 surfaces and abelian surfaces the numbers can 
be expressed solely in terms of quasimodular forms
might be related to physical dualities.   
The numbers $t^S_\delta(L)$ will on a general surface (including also 
$\P_2$) count all curves with the prescribed numbers of nodes, including the
reducible ones. It seems however that on an abelian surface one 
is actually counting irreducible curves, and in the case of K3 
surfaces one can restrict attention to the case that
$|L|$ contains only irreducible reduced curves.
Both in the case of abelian surfaces and K3 surfaces I do then expect the
generating function to count the curves of given genus in sub-linear systems
of $|L|$, even if $L$ is only assumed to be very ample and not a high multiple
of an ample line bundle. The curves can then have worse singularities than
nodes and a curve $C$ should be counted  with a multiplicity determined 
by the singularities of $C$, as in \cite{B}, \cite{F-G-vS}.
In particular nodal (or more generally immersed) curves should count with
 multiplicity $1$.  In the case of K3 surfaces and primitive
line bundle $L$ the numbers of these curves have  in the meantime been
computed in \cite{Br-Le}.

The coefficients of the unknown power
series are by the recursion of \cite{C-H} the solutions to a highly 
overdetermined system of linear equations, the same is true for a similar 
recursion  obtained by Vakil \cite{Va} for rational ruled surfaces and 
the results of \cite{Ch} on $\P_2$ and $\P_1\times \P_1$. This gives an 
additional check of the conjecture.
Finally we compute the numbers of genus $2$ curves on an abelian surface.

I thank B. Fantechi for many useful discussions without which this  
paper could not have been written and P. Aluffi for pointing out \cite{V2} 
to me. I thank Don Zagier for useful comments and discussions which improved 
the formulation of conjecture \ref{mainconj}, 
R. Vakil for giving me a preliminary version of \cite{Va} and
S. Kleiman, R. Piene  for useful comments. 
This paper was started during my stay at the  Mittag Leffler Institute.

\section{Statement of the conjecture}

Let $S$ be a projective algebraic surface and $L$ a line bundle on $S$.  
In this paper by a curve on $S$ we mean an effective reduced divisor on $S$. 
A nodal curve on $S$ is a reduced (not necessarily irreducible) divisor on $S$,
which has only nodes as singularities.
We denote by $K_S$ the  canonical bundle and by $c_2(S)$ the degree
of the second Chern class. For two line bundles $L$ and $M$ let
$LM$ denote the degree of $c_1(L)\cdot  c_1(M)\in H^4(S,\Z)$.  

In \cite{V2} (for $\delta\le 6$) and \cite{K-P} (for $\delta\le 8$) formulas
for the  number $a^S_{\delta}(L)$ of $\delta$-nodal curves  in a general 
$\delta$-dimensional linear sub-system $V$ of $|L|$ were proved. 
Here general means that $V$ lies in an open subset of the Grassmannian of 
$\delta$-dimensional subspaces of $|L|$. 
The number is expressed as a  polynomial in
$c_2(S)$, $K_S^2$, $LK_S$ and  $L^2$ of degree $\delta$. 
The formulas are valid if $L$ is a sufficiently high 
multiple of an ample line bundle. In other words, for such $L$, 
the locally closed subset $W^S_\delta(L)$ (with the reduced structure)  
of elements in $|L|$ defining  $\delta$-nodal curves  has codimension
$\delta$, and its degree is given by a polynomial as above.

It is clear and was noted before  \cite{K-P}
that there should be a formula for all $\delta$:

\begin{conj} \label{aconj}
For all $\delta\in  \Z$ there exist universal polynomials
$T_\delta(x,y,z,t)$ of degree $\delta$ ($T_\delta=0 $ 
if $\delta<0$) with the following property. Given $\delta$  and a pair 
of a surface $S$ and a very ample line bundle $L_0$ there exists an 
$m_0>0$ such that for all $m\ge m_0$  and for all  very ample line 
bundles $M$ the line bundle $L:=L_0^{\otimes m}\otimes M$ satisfies
$$
a^S_{\delta}(L)=T_\delta(L^2,LK_S,K_S^2,c_2(S) ).
$$
\end{conj} 

\begin{rem} Note that the statement is slightly stronger than that of 
\cite{V2} in the case $\delta\le 6$. We expect $L$ does not have to be 
  a high power of a very ample line bundle but that it suffices that 
$L$ is sufficiently ample. In fact the results of the final two sections below
suggest that there might exist a universal constant $C >0$ 
(independent of $S$ and $L$), such that, if $L$ is $C\delta$-very ample 
(see section 5), then the conjecture holds for  up to $\delta$-nodes.
\end{rem}

Assuming  conjecture \ref{aconj} we will in future just write
$$
T^S_\delta(x,y):=T_\delta(x,y,K_S^2,c_2(S) ),\quad  
t^S_\delta(L):=T^S_\delta(L^2,LK_S)),
\quad T(S,L):=\sum_{\delta\ge 0} t^S_\delta(L) x^\delta.$$
The aim of this note is to give a conjectural formula for the 
generating function $T(S,L)$, and to give some evidence for it.

We start by noting that conjecture 2.1 imposes rather strong 
restrictions on the structure of
$T(S,L)$. The point is that the conjecture applies to all surfaces, 
including those with several connected components. In  this case we will write
$|L|$ for $\P(H^0(L))$. By definition  $W^S_\delta(L)$  includes only
$f\in |L|$ which do not vanish identically on a component of $S$.

\begin{prop}\label{constr}
Assume conjecture \ref{aconj}. Then there exist universal power series
$A_1$, $A_2$, $A_3$, $A_4\in \Q[[x]]$ such that 
$$T(S,L)=\exp(L^2 A_1+LK_S A_2+ K_S^2 A_3 +c_2(S) A_4).$$
\end{prop}
\begin{pf}
Fix $\delta_0\in \Z_{>0}$. Assume that $S=S_1\sqcup S_2$ and that 
$L_1:=L|_{S_1}$ and $L_2:=L|_{S_2}$  are both sufficiently ample so that the
$W^{S_i}_\delta(L_i)$ have codimension $\delta$ and degree 
$t^{S_i}_\delta(L_i)$ in $|L_i|$ for $i=1,2$ and $\delta<\delta_0$.
Fix $\delta<\delta_0$.
The application $(f+g)\mapsto (f,g)$ defines a surjective morphism  
$p:U\to |L_1|\times |L_2|$, defined on the open subset $U\subset |L|$ where 
neither  $f$ not $g$ vanish identically. The fibres of $p$ are lines in $|L|$.
Obviously 
$$
W^S_\delta(L)=p^{-1}\Big(\coprod_{\delta_1+\delta_2=\delta}
W^{S_1}_{\delta_1}(L_1)\times W^{S_2}_{\delta_2}(L_2)\Big).
$$
In particular $W^S_\delta(L)$ has codimension $\delta$ in $|L|$
and modulo the ideal generated by $x^{\delta_0}$ we have
$T(S,L)\equiv T(S_1,L_1) T(S_2,L_2).$

Now choose $n\in \Z_{>0}$ such that conjecture \ref{aconj} 
holds for $Z_1:=(\P_2,\oo(n))$, $Z_2:=(\P_2,\oo(2n))$, 
$Z_3:=(\P_2,\oo(3n))$ and
$Z_3:=(\P_1\times
\P_1,\oo(n,n))$ for all $\delta\le \delta_0$.
Let $S=S(a_1,a_2,a_3,a_4)$ be the disjoint union of $a_i$ copies of each of 
the $Z_i$ with
$a_i\in Z_{\ge 0}$. Then by the above argument $T(S,L)=\prod_i T(Z_i)^{a_i}$ 
which  can be
written in the form $\exp(L^2 A_1+LK_S A_2+ K_S^2 A_3 +c_2(S) A_4),$
for universal power series $A_1$, $A_2$, $A_3$, $A_4$.
On the other hand we note that the $4$-tuple
$(L^2,LK_S,K_S^2,c_2(S))$ takes on the $Z_i$ the linearly independent
values $(n^2,-3n,9,3)$, $(4n^2,-6n,9,3)$, $(9n^2,-9n,9,3)$, $(2n^2,-4n,8,4)$.
Thus the image of  the $S(a_1,a_2,a_3,a_4)$  is Zariski-dense   in
$\Q^4$, and the result follows.
\end{pf}

 We   recall some facts about quasimodular forms from
\cite{K-Z}. We denote by    $\H:=\bigl\{\tau\in \C\bigm|\Im(\tau)>0\bigr\}$   
the  complex upper half plane, and for $\tau \in \H$ we write 
$q:=e^{2\pi i\tau}$.
A modular form of weight $k$ for $SL(2,\Z)$ is a holomorphic function 
$f$ on $\H$ satisfying
$$
f\Big(\frac{a\tau + b}{c\tau + d}\Big)=(c\tau+d)^k f(\tau), 
\quad \tau\in \H,\ 
\quad\left(\begin{matrix}a&b\\c&d\end{matrix}\right)\in SL(2,\Z)
$$
and having a  Fourier series
$f(\tau)=\sum_{n=0}^\infty a_n q^n.$

Writing $\sigma_k(n):=\sum_{d|n}d^{k}$, the Eisenstein series 
$$
G_k(\tau)=-\frac{B_k}{2k}+\sum_{n>0}\sigma_{k-1}(n)q^n,\quad
k\ge 2, \quad B_k=\text{$k$th Bernoulli number}
$$
are for even $k>2$ modular forms of weight $k/2$,
while $G_2(\tau)$ is only a quasimodular form. 
Another important modular form is the discriminant
$$
\Delta(\tau)=q\prod_{k>0} (1-q^k)^{24}=\eta(\tau)^{24}
$$ 
where $\eta(\tau)$ is the Dedekind $\eta$ function.

For the precise definition of quasimodular forms see  \cite{K-Z}. 
They are essentially the holomorphic parts of almost holomorphic
 modular forms. 
The ring of  modular forms for $SL(2,\Z)$ is just 
$\Q[G_4,G_6]$, while the ring of quasimodular forms 
is $\Q[G_2,G_4,G_6]$. 
We denote by $D$ the differential operator 
  $D:=\frac{1}{2\pi i}\frac{d}{ d\tau}=q\frac{d}{d\,q}$.
Unlike the ring  of  modular forms the ring of quasimodular forms is
 closed under differentiation, i.e. for a quasimodular form $f$ of 
weight $k$ the derivative 
$Df$ is a quasimodular form of weight $k+2$. In fact every quasimodular
form has a unique representation as a sum of 
derivatives of modular forms and of $G_2$ (see \cite{K-Z}).

\begin{conj}\label{mainconj} There exist universal power series $B_1$, $B_2$ 
in $q$ such that 
$$
\sum_{\delta\in \Z}t^S_\delta(L) (DG_2(\tau))^{\delta}=
 \frac{(DG_2(\tau)/q)^{\chi(L)}\,B_1(q)^{K_S^2}\, B_2(q)^{LK_S}}
{(\Delta(\tau)\,D^2G_2(\tau)/q^2)^{\chi(\oo_S)/2}}.
$$
\end{conj}

\begin{rem}
\noindent (1) 
Using the fact that $DG2(\tau)/q$, $B_1(q)$, 
$B_2(q)$, $\Delta(\tau)D^2G_2(\tau)/q^2$ are power series in $q$ starting
with $1$, and by the standard formulas 
  $\chi(\oo_S)=(K_S^2+c_2(S))/12$, $\chi(L)=(L^2-LK_S)/2+\chi(\oo_S)$
one sees that conjecture \ref{mainconj} 
expresses the $t^S_\delta(L)$ as polynomials of degree $\delta$ in 
$L^2,\, K_SL,\, K_S^2,\, c_2(S)$.

\noindent (2)
I have checked that conjecture \ref{mainconj} reproduces the formulas of
Vainsencher and Kleiman-Piene for $\delta\le 8$. This determines $B_1(q)$ and 
$B_2(q)$ up to degree  $q^8$. In  remark \ref{Srem} below we   use the
  formulas of \cite{C-H} for the Severi degrees in $\P_2$ to determine
 the coefficients of $B_1(q)$ and 
$B_2(q)$ up to degree 28 (they are given here  up to degree 20). 

\begin{align*}   B_1&(q)\, \equiv\, 
1-q-5\,{q}^{2}+39\,{q}^{3}-345\,{q}^{4}+2961\,{q}^{5}-24866\,{q}^{6}+
207759\,{q}^{7}-1737670\,{q}^{8}\\ &+14584625\,{q}^{9}-122937305\,{q}^{10}
+1040906771\,{q}^{11}-8852158628\,{q}^{12}+75598131215\,{q}^{13}\\ &-
648168748072\,{q}^{14}+5577807139921\,{q}^{15}-48163964723088\,{q}^{16
}+417210529188188\,{q}^{17}\\ &-3624610235789053\,{q}^{18}+
31575290280786530\,{q}^{19}-275758194822813754\,{q}^{20}+O  ({q}^{ 21}  )\\
B_2&(q)\, \equiv \, 
1+5\,q+2\,{q}^{2}+35\,{q}^{3}-140\,{q}^{4}+986\,{q}^{5}-6643\,{q}^{6}+
48248\,{q}^{7}-362700\,{q}^{8}\\ &+2802510\,{q}^{9}-22098991\,{q}^{10}+
177116726\,{q}^{11}-1438544962\,{q}^{12}+11814206036\,{q}^{13}\\ &-
97940651274\,{q}^{14}+818498739637\,{q}^{15}-6888195294592\,{q}^{16}+
58324130994782\,{q}^{17}\\ &-496519067059432\,{q}^{18}+4247266246317414\,{
q}^{19}-36488059346439524\,{q}^{20}+O(q^{21}).
\end{align*}
\end{rem}

\begin{rem}
We give a reformulation of the conjecture.
We define, for all $l,m,r\in \Z$  
\begin{align*}
n^S_r(l,m)&:=T^S_{l+\chi(\oo_S)-1-r}(2l+m,m)\\
m^S_g(l,m)&:=n_{g-m-2+\chi(\oo_S)}(l,m)
\end{align*}
If $L$ is sufficiently ample with respect to $\delta=\chi(L)-1-r$ and $S$
(and thus in particular $\chi(L)=H^0(S,L)$ and $r \ge 0$), then 
$n^S_{r}((L^2-LK_S)/2,LK_S)$ counts  the 
$\delta$-nodal curves in a general $r$-codimensional sub-linear system of 
$|L|$. Then
$$
\sum_{l\in \Z} n_r^S(l,m)q^{l} =B_1(q)^{K_S^2}B_2(q)^{m}\big(DG_2(\tau)\big)^r
\frac{D^2 G_2(\tau)}{(\Delta(\tau)D^2G_2(\tau))^{\chi(\oo_S)/2}}.
$$
An irreducible  $\delta$-nodal curve $C$ on  $S$ has geometric
genus $g(C)=(L^2+LK_S)/2+1-\delta$, we take the same definition
also if $C$ is reducible. For $L$ sufficiently ample
$m^S_g((L^2-LK_S)/2,LK_S)$ counts  the nodal curves $C$  with $g(C)=g$ in a
general $g-LK_S  +\chi(\oo_S)-2$-codimensional sub-linear system of $|L|$.
\end{rem}
\begin{pf}
If $f(q)$ and $g(q)$
are power series in $q$ and $g(q)$ starts with $q$, then we can develop
$f(q)$ as a power series in $g(q)$ and 
$$
\text{Coeff}_{g(q)^k}
f(q)=\text{Res}_{g(q)=0}\frac{f(q)dg(q)}{g(q)^{k+1}}
=\text{Coeff}_{q^{0}}\frac{f(q)Dg(q)}{g(q)^{k+1}}.$$
We apply this with $g(q)=DG_2(\tau)$.
\end{pf}

\section{Counting curves on K3 surfaces and abelian surfaces}

Let now $S$ be a surface with numerically trivial canonical divisor.
We denote $n^S_r(l):=n^S_r(l,0)$, i.e. for $L$ sufficiently ample
$n^S_r(L^2/2)$ is the number
of $\chi(L)-r-1$-nodal curves in an $r$-codimensional sub-linear system of $L$.
$n^S_{r}(l)$ can be expressed in terms of quasimodular forms. 

For $S$ a K3 surface, $A$ an abelian surface and $F$ an Enriques or 
bielliptic surface we get 
\begin{align} \label{fK3}
\sum_{l\in\Z} n^S_{r}(l)q^l &= \big(DG_2(\tau)\big)^r/\Delta(\tau)\\
\label{fA}  
\sum_{l\in \Z} n^A_{r}(l)q^l&= \big(DG_2(\tau)\big)^r D^2G_2(\tau),\\
\sum_{l\in \Z,r\ge 0} n^A_{r}(l)q^l\frac{z^r}{r!}
&=\frac{1}{z}D\big(\exp\big(D G_2(\tau) z\big)\big)\nonumber \\
\sum_{l\in \Z,r\ge 0} n^F_{r}(l)q^l &=\big(DG_2(\tau)\big)^r
\big( (D^2G_2(\tau))/\Delta(\tau)\big)^{1/2}
\end{align}
Note that 
$m_g^S(l)=n_g^S(l), \ m_g^F(l)=n_{g-1}^F(l), \ m_g^A(l)=n_{g-2}^A(l).$ 

\begin{rem}\label{exprem}  In the case of an abelian surface or a K3 surface 
we expect that the numbers $n^A_{r}(l)$ and  $n^S_{r}(l)$ have a
more interesting geometric significance.

\noindent(1) 
Let $(S,L)$ be a polarized K3 surface with 
$Pic(S)=\Z L$. Then the linear system $|L|$ contains 
only irreducible curves. The number  $n^S_r(L^2/2)$ is a count of 
curves $C\in |L|$ of geometric genus $r$ passing through $r$ general points on 
$S$.
In this case the numbers of rational curves have been calculated in \cite{Y-Z} 
and \cite{B} and (\ref{fK3}) is  a generalization to arbitrary genus.
The rational curves that are counted  are not necessarily
  nodal. In the count a rational  curve $C$ is assigned 
the Euler number $e(\bar JC)$ of its compactified Jacobian 
(which is $1$ if $C$ is 
immersed) as   multiplicity. 

Let $\overline M_{g,n}(S,\beta)$ be the moduli
space of $n$-pointed genus $g$ stable maps of homology class $\beta$
(see \cite{K-M}).
It comes equipped with an evaluation map $\mu$ to $S^n$.
In \cite{F-G-vS} it is shown that for a rational curve $C$ on a K3 surface 
$e(\bar JC)$ is just the multiplicity of $\overline M_{0,0}(S,[L])$ at 
the point  defined by the  normalization of $C$. Here  $[L]$ denotes the 
homology class Poincar\'e dual to $c_1(L)$.
In other words $n^S_{0}(L^2/2)$ is just the length of the 
$0$-dimensional scheme $\overline M_{0,0}(S,[L])$.
I expect that for curves of arbitrary genus the corresponding result
should hold: If $S$ and $L$ are general and $x$ is a general point in 
$S^r$ then the fiber $\mu^{-1}(x)\subset \overline M_{r,r}(S,[L])$ 
should be a finite scheme and $n^S_{r}(L^2/2)$ should just be its length. 
More generally  
 $n_r^S(L^2/2)$ should be  a generalized Gromov-Witten invariant as 
defined and
studied  in \cite{Br-Le} in the symplectic setting and in \cite{Be-F2} in the
algebraic geometric setting. In the meantime this invariant 
has been computed in \cite{Br-Le}  for curves of arbitrary genus on K3 
surfaces
confirming the conjecture.

\noindent(2) 
Let $A$ be an  abelian surface with a very  ample line bundle $L$. 
We claim that in general all the curves counted in $n^A_{r}(L^2/2)$ 
will be irreducible and reduced:
The set of $\delta$-nodal curves in $|L|$ has expected dimension 
$\chi(L)-\delta-1=L^2/2-\delta-1$. On the other hand the set of 
reducible $\delta$-nodal curves $C_1+\ldots +C_n\in |L|$ with 
$C_i\in |L_i|$ has expected dimension
$$
\sum_{i=1}^n(L_i^2/2 -1 ) -\delta+\sum_{1\le i\ne j\le n}L_i L_j=
L^2/2-\delta-n.
$$
Therefore in general $n^A_r(L^2/2)$ should  count the irreducible curves 
$C\in |L|$  of geometric genus $r+2$ passing through $r$ general points.
Again we expect that this result holds in a modified form also if $L$ is 
not required to be sufficiently   ample, and if not 
all the curves in $|L|$ are immersed.
The moduli space $\overline M_{g,n}(A,[L])$ is naturally fibered over 
$Pic^L(A)$. The fibers are the spaces  $\overline M_{g,n}(A,|M|)$
of stable maps $\phi: W\to A$ with $\phi_*(W)$ a divisor in  $|M|$, where
$c_1(M)=c_1(L)$. Again for $A$ and $L$ general and  $x$ a general point in 
$A^r$ the
number 
$n^A_{r}(L^2/2)$ should be the length of the fiber 
$\mu^{-1}(x)\subset \overline M_{r+2,r}(A,|L|)$.
\end{rem}

We want to  show conjecture \ref{mainconj} and the expectations of  remark
\ref{exprem} for   abelian surfaces in a special case.

\begin{thm} 
Let $A$ be an   abelian surface with an ample line bundle $L$ such that 
$c_1(L)$ is a
polarization of type
$(1,n)$.   Assume that
$A$ does not contain elliptic curves. Write again $\sigma_1(n)=\sum_{d|n} d$.
Then the number of  genus
$2$ curves in
$|L|$ is $n^2\sigma_1(n)$. Moreover all these curves are irreducible and 
immersed, and the moduli space
$\overline M_{2,0}(A,|L|)$ consists of $n^2\sigma_1(n)$ points
corresponding to the their normalizations.
\end{thm}
\begin{pf}
We will denote by $[C]$ the homology class of a curve $C$ and
by 
$[L]$ the Poincar\'e dual of $c_1(L)$. For a divisor $D$ we write $c_1(D)$ for 
$c_1(\oo(D))$.
Let $\phi:C\to A$ be a morphism from a connected 
nodal curve of arithmetic genus $2$ to 
$A$, with $\phi_*([C])=[L]$.
Put $D:=\phi(C)$. As $A$ does not contain curves of genus $0$ or $1$,
$C$ must be  irreducible and smooth, and $\phi$ must be generically injective.
In particular $[D]=[L]$.
Let $J(C)$ be the Jacobian of $C$. We freely use  standard results about
Jacobians of curves, see (\cite{L-B} chap. 11) for reference.
For each $c\in C$   the Abel-Jacobi map $\alpha_c:C\to J(C)$ is an embedding
with $\alpha_c(c)=0$.
We write $C_c:=\alpha_c(C)$ and $\theta_C:=c_1(\alpha_c(C))$.
By the Torelli theorem the isomorphism class of $C$ is determined uniquely by
the pair $(J(C),\theta_C)$. For all $a\in A$ we denote by $t_a$ the
translation  by $a$. By the universal property of the Jacobian there is  a 
unique isogeny
$\widetilde \phi:J(C)\to A$ such that 
$\phi=t_{\phi(c)}\circ  \widetilde\phi \circ \alpha_c$ for all $c\in C$.
As $\widetilde \phi$ is \'etale $\phi$ is an immersion and $\phi:C\to D$ is the
normalization map.
We also see that $\widetilde\phi^*(c_1(L))=n\theta_C$.

On the other hand, let $(B,\gamma)$ be a principally polarized abelian 
surface and  $\psi:B\to A$ an isogeny with $\psi^*(c_1(L))=n\gamma$.
By the criterion of Matsusaka-Ran and the   assumption that $A$ does not
contain elliptic curves we obtain that $(B,\gamma)=(J(C),\theta_C)$ for
$C$ a smooth curve of genus $2$ and $\psi=\widetilde \phi$ for   a morphism
$\phi:C\to A$  with $\phi_*([C])=[L]$.
$\widetilde \phi$ depends only on $\phi$ up to composition with a 
translation in $A$ and 
$\phi=t_{\phi(c)}\circ\widetilde\phi\circ\alpha_c$ is determined by 
$\widetilde\phi$ up to translation in $A$.
 
By  the universal property of $J(C)$ applied to the embedding $\alpha_c$, an
automorphism $\psi$ of
$C$ induces an automorphism $\widehat \psi$ of
$J(C)$. 
If $\epsilon$ is an automorphism of $J(C)$ with
$\epsilon^*(\theta_C)=\theta_C$, then it is
$\widehat\psi$ for some automorphism $\psi$ of $C$.
($H^0(J(C),\oo(C_c))=1$, and  $a\mapsto t_a^*(\oo(C_c))$ 
defines an isomorphism  $J(C)\to Pic^{\theta_C}(J(C))$.)

Therefore we see that  the set $M_1$ of morphisms $\phi:C\to A$ 
from curves of genus $2$ with
$\phi_*([C])=[L]$ modulo composition with  automorphisms of $C$ and with  
translations of $A$ can be identified with the
set $M_2$ of morphisms $\psi:B\to A$ from a principally polarized 
abelian surface $(B,\gamma)$, such that $\psi^*(c_1(L))=n\gamma$ modulo 
composition with automorphisms 
$\eta:B\to B$ with $\eta^*(\gamma)=\gamma$.

We write $A=\C^2/\Gamma$ and $B=\C^2/\Lambda$.
Then $c_1(L)$ is given by an alternating form $a:\Gamma\times\Gamma\to \Z$ 
such that there is a basis $x_1,x_2,y_1,y_2$ of $\Gamma$ with 
$a(x_1,y_1)=1$, $a(x_2,y_2)=n$, 
$a(x_1,x_2)=a(y_1,y_2)=0$.
A homomorphism $\psi:B\to A$ is given by a linear map 
$\widehat \psi:\C^2\to \C^2$
with $\psi(\Lambda)\subset \Gamma$. We see that 
$M_2$ can be identified with the set $M_3$ of sublattices 
$\Lambda\subset \Gamma$
of index $n$ with $a(\Lambda,\Lambda)\subset n\Z$.
We claim that $M_3$ has $\sigma_1(n)$ elements. 

First we want to see that this claim implies the theorem. 
Let $Pic^L(A)$ be the group of line bundles on $A$ with first Chern class 
$c_1(L)$.
By \cite{L-B} proposition 4.9 the morphism
$\phi_L:A\to Pic^L(A); a\mapsto t^*_aL$ is \'etale of degree $n^2$.
By the claim this means that for each $L_1\in Pic^L(A)$ the linear
system $|L_1|$ contains precisely $n^2\sigma_1(n)$ curves of genus $2$.

Finally we show the claim.
Via the basis $x_1,y_1,x_2,y_2$ (in that order) we identify $\Gamma$ with
$\Z^2\times \Z^2$. We see that $\Lambda$ must be of the form 
$\Lambda'\times \Z^2$, where $\Lambda'$ is a sublattice of index $n$ in 
$\Z^2$, satisfying $b(\Lambda',\Lambda')\subset n\Z$, for the alternating form  
$b$  defined by $a(x_1,y_1)=1$. Let $\Lambda'$ be a sublattice of $\Z^2$
of index $n$. We claim that  $b(\Lambda',\Lambda')\subset
n\Z$. Then the result follows by the well-known fact that the number of
sublattices of index $n$ in a rank two lattice is $\sigma_1(n)$.
Let $L_1:=p_2(\Lambda')$ for the second projection
$\Z^2\to \Z$ and put $L_2:=ker(p_2|_{\Lambda'})$.
Then $L_1=d_1\Z$ and $L_2=d_2\Z$ for $d_1,d_2\in \Z$ with $d_1d_2=n$.
Choose $x=(k,d_1)\in \Lambda'\cap p_2^{-1}(d_1)$. Then 
$\Lambda'$ is generated by $L_2$ and $x$, and in particular
$b(\Lambda',\Lambda')\subset n\Z$
\end{pf}

\section{Severi degrees on $\P_2$ and rational ruled surfaces}

The Severi   degree
$N^{d,\delta}$ is the number of plane curves of degree $d$
with $\delta$ nodes passing through $(d^2+3d)/2-\delta$ general 
points. In \cite{R2} a recursive procedure for determining the
$N^{d,\delta}$ is shown, and in \cite{C-H} a different recursion
formula is proven. 
In the number $N^{d,\delta}$ also reducible curves are included,
they  however occur only if $d\le \delta+1$, furthermore   the numbers of
irreducible curves can be determined from them (\cite{C-H} see also \cite{Ge}).
For simplicity I will write $t_\delta(d):=t^{\P_2}_{\delta}(\oo(d))$.
If $d$ is sufficiently large with respect to $\delta$, then 
$N^{d,\delta}$ should be equal to $t_{\delta}(d)$:
 
\begin{conj}\label{Sconj}
 If $\delta\le 2d-2$, then 
$N^{d,\delta}=t_{\delta}(d)$.
\end{conj}

\begin{rem}\label{Srem}

(1) 
Conjecture 
\ref{Sconj} and \cite{C-H}  provide  an effective
 method of determining the coefficients of the two unknown power series 
$B_1(q)$ and $B_2(q)$.
Using  a suitable program I computed  the $N^{d,\delta}$ via the 
recursive  formula from   \cite{C-H} for $d\le 16$ and $\delta\le 30$.
 We write $x=DG_2(\tau)$. By conjectures \ref{mainconj} 
and \ref{Sconj} one has for all $d>0$ moduli the ideal generated by
$x^{2d-1}$ the identity
\begin{equation}\label{expeqn}
\sum_{\delta\in \Z}N^{d,\delta}x^\delta\equiv \exp(d^2C_1(x)+dC_2(x)+C_3(x))
\end{equation}
Here $C_1(x)$ is known by conjecture \ref{mainconj} and the first $k$ 
coefficients of $C_2(x) $ and $C_3(x)$ determine the first $k$ coefficients
of $B_1(q)$ and $B_2(q)$. Taking logarithms on both sides gives, for 
any two degrees
$d_1<d_2$ and $\delta\le 2d_1-2$, a system of  two linear equations for the 
coefficients of $x^\delta$   in $C_2(x)$ and $C_3(x)$.   
Note that this also gives a
test of the conjecture. It is a nontrivial fact that the generating function
has the special shape (\ref{expeqn}). In particular each pair
$d_1<d_2$ with $\delta\le 2d_1-2$ already determines the coefficients
of $x^\delta$.

\noindent(2) The conjecture implies in particular that for 
$\delta\le 2d-2$ the numbers $N^{d,\delta}$ are given by a
polynomial
$P_\delta(d)$ of
degree $2\delta$ in $d$. This was already conjectured in \cite{D-I}. Denote by
$p_\mu(\delta)$ the coefficient of
$d^{2\delta-\mu}$ in $P_\delta$.
 In 
\cite{D-I} a conjectural formula for the leading coefficients
$p_\mu(\delta)$ for $\mu\le 6$ is given. 
Kleiman and Piene have determined $p_7(\delta)$ and $p_8(\delta)$ 
\cite{K-P}. In \cite{Ch} 
the  $N^{d,\delta}$ for $\delta\le 4$ are computed using the recursive method 
of \cite{R2}, and as an application $p_0( \delta)$ and $p_1(\delta)$
are determined.  Using conjectures \ref{mainconj} and  \ref{Sconj} there is an
algorithm to determine  the $p_\mu(\delta)$. Again
we use the formula
(\ref{expeqn}), and collect terms. From knowing the coefficients of $C_2(x)$
and $C_3(x)$ up to degree 28 we get
 the $p_\mu(\delta)$ for $\mu\le 28$. 
For $\mu\le 8$ they coincide with those from
\cite{D-I},\cite{Ch} and
\cite{K-P}.  Let $[\ ]$ denote the integer part. 
For $\mu\le 28$ we observe that $p_\mu(\delta)$ is of the form
$$p_\mu(\delta)=\frac{3^{\delta-[\mu/2]}}{(\delta-[\mu/2])! }Q_\mu(\delta),$$
where $Q_\mu(\delta)$ is a polynomial of degree $[\mu/2]$ 
in $\delta$ with integer coefficients, which have only products of 
powers of $2$ and $3$ as common factors. In particular 
\begin{align*}
Q_8(\delta)&=-2^4(282855\, \delta^4-931146\, \delta^3
+417490\, \delta^2+425202\, \delta+1141616),\\
Q_9(\delta)&=-2^33^2(128676\, \delta^4+268644\, \delta^3-1011772\, \delta^2-
1488377\, \delta-1724779),\\
Q_{10}(\delta)&=2^4
3^2(4345998\, \delta^5-15710500\, \delta^4-3710865\, \delta^3+7300210\, 
\delta^2\\
&\qquad +57779307\, \delta+98802690).
\end{align*}
Note that for $d>\delta+1$ all the curves are irreducible,
so that in this case we get    a conjectural formula for the 
stable Gromov-Witten invariants of $\P_2$, i.e. the numbers of
irreducible curves  $C$ of degree $d$ with
$g(C)\ge\binom{d-2}{2}$.
\end{rem}

\begin{rem}\label{sigmarem} Let $\Sigma_e$ be a rational ruled surface, 
$E$ the curve with
$E^2=-e$ and $F$ a fiber of the ruling. For simplicity denote
$t^e_\delta(n,m):=t^{\Sigma_e}_{\delta}(\oo(nF+mE))$. \cite{Ch} determined the
$N^{2,m}_{0,\delta}$ for $\delta\le 9$ as polynomials of degree $\delta$ in $m$
and several numbers $N^{3,m}_{0,\delta}$
In \cite{Va} a
recursion formula  very similar to that
of \cite{C-H} is proved for the generalized Severi degrees 
$N^{n,m}_{e,\delta}$, i.e. the number 
of $\delta$-nodal curves in $|nF+mE|$ which do not contain $E$ as 
a component. Using a suitable program I computed the $N^{n,m}_{e,\delta}$
for $e\le 4 $, $\delta\le 10$, $n\le 11$ and $m\le 8$.
These results are compatible with conjecture \ref{mainconj}, if one conjectures
that for $(n,m)\ne (1,0)$ one has
$N^{n,m}_{e,\delta}=t^{e}_{\delta}(n,m)$ if and only if 
$\delta\le min(2m,n-em)$ or $\delta\le min(2m,2n)$ in case 
$e=0$. 
\end{rem}

\begin{rem}
A slightly sharpened version of conjecture \ref{Sconj} can be reformulated
as saying that 
$N^{d,\delta}=t_{\delta}(d)$ if and only if $H^0(\P_2,\oo(d))>\delta$ and 
the locus of nonreduced curves in $|\oo(d)|$  has codimension bigger then
$\delta$.
In a similar way  the conjecture of remark \ref{sigmarem}
 can be reformulated as
saying that $N^{n,m}_{e,\delta}=t^{e}_{\delta}(n,m)$ if and only if
$H^0(\Sigma_e,\oo(nF+mE))>\delta$ and 
the locus of curves in $|\oo(nF+mE)|$ which are nonreduced or contain $E$ as a
component has codimension bigger then $\delta$.
One would expect that nonreduced curves   contribute to 
the count of nodal curves, and the recursion formula of \cite{Va}
only counts curves not containing $E$. Therefore it seems that, at 
least in the case of
$\P_2$ and of rational ruled surfaces, $t^{S}_{\delta}(L)$ is the actual
number of $\delta$-nodal curves in a general $\delta$-codimensional linear
system, unless this cannot be expected for obvious   geometrical reasons.
\end{rem}

\section{Connection with Hilbert schemes of points}

Let again $S$ be an algebraic surface and let $L$ be a line bundle on $S$.
Let $S^{[n]}$ be the Hilbert scheme of finite subschemes of length $n$ on $S$,
and let $Z_n(S)\subset S\times S^{[n]}$ denote the universal family with 
projections 
$p_n:Z_n(S)\to S$, $q_n:Z_n(S)\to S^{[n]}$. 
Then
$L_n:=(q_n)_*(p_n)^*(L)$ is a locally free sheaf of rank $n$ on $S^{[n]}$.

\begin{defn}
Let $S^{\delta}_2\subset S^{[3\delta]}$ be the closure (with the reduced
induced structure) of the locally closed subset $S^{\delta}_{2,0}$ which 
parametrizes
subschemes of the form
$\coprod_{i=1}^{\delta} \spec(\oo_{S,x_i}/\mmm^2_{S,x_i})$, where $x_1,\ldots
x_{\delta}$ are distinct points in
$S$. It is easy to see that $S^{\delta}_2$ is birational to $S^{[\delta]}$.
We put 
$
d_n(L):=\int_{S^\delta_2}c_{2\delta}(L_{3\delta}).
$
\end{defn}

Following \cite{B-S} we call $L$
$k$-very ample if for all subschemes $Z\subset S$ of length $k+1$ the natural
map
$H^0(S,L)\to H^0(L\otimes\oo_Z)$ is surjective. If $L$ and $M$ are very ample, 
then
$L^{\otimes k}\otimes M^{\otimes l}$ is
$(k+l)$-very ample.

\begin{prop} \label{port} Assume $L$ is $(3\delta-1)$-very ample, then   a 
general  $\delta$-dimensional linear subsystem $V\subset |L|$ contains  
only finitely many curves
$C_1,\ldots,C_s$ with   $\ge \delta$ singularities.
There   exist positive integers
$n_1,\ldots,n_s$ such that
$\sum_i n_i=d_\delta(L)$.  
If furthermore $L$ is $(5\delta-1)$-very ample ($5$-very ample if $\delta=1$), 
then  the
$C_i$ have precisely
$\delta$ nodes as singularities.
\end{prop}
\begin{pf} 
Assume first that $L$ is $(3\delta-1)$-very ample.
We apply the Thom-Porteous formula  to the
restrictions of the  evaluation map
$H^0(S,L)\otimes\oo_{S^{[3\delta]}}\to L_{3\delta}$ to $S^{\delta}_2$ and to
$S^{\delta}_{2}\setminus S^{\delta}_{2,0}$. As $L$ is $(3\delta-1)$-very ample 
the evaluation map is surjective. Then (\cite{Fu}
ex. 14.3.2) applied to
$S^{\delta}_2$ gives that for a general $\delta$-dimensional sublinear system 
$V\subset |L|$ the class $d_n(L)$ is represented by the class  of the 
finite scheme $W$ of 
$Z\in S^{\delta}_2$ with $Z\subset D$ for $D\in V$. The scheme structure of 
$W$ might
be nonreduced.  The application of (\cite{Fu} ex. 14.3.2) to
$S^{\delta}_{2}\setminus S^{\delta}_{2,0}$  and a dimension count give that 
$W$ lies entirely in
$S^{\delta}_{2,0}$.

Now assume that $L$ is $(5\delta-1)$-very ample. 
Let $V\subset |L|$ again be general $\delta$-dimensional subsystem of $|L|$.
The Porteous formula applied to the restriction of $L_{3\delta+3}$ to
$S^{\delta+1}_2$ and a  dimension count shows that
there will be no curves in $V$ with more than $\delta$ singularities.

Let $S^{\delta}_{3,0}\subset S^{[5\delta]}$ be the locus of 
schemes of the form $Z_1\sqcup Z_2\ldots \sqcup Z_{\delta}$, where each 
$Z_i$ is of the form $\spec(\oo_{S,x_i}/(\mmm^3+xy))$ with $x,y$
local parameters at $x_i$ and let $S^{\delta}_{3}$ be the closure. 
If a curve $C$ with precisely $\delta$ singularities 
does not contain a subscheme
corresponding to a point in
$S^{\delta}_{3}\setminus S^{\delta}_{3,0}$,
then it has
$\delta$ nodes as only singularities. It is easy to see that
$S^{\delta}_{3,0}$ is smooth of dimension $4\delta$. 
Applying   the Porteous formula to the restriction of $L_{5n}$
to  $S^{\delta}_{3}\setminus S^{\delta}_{3,0}$ and a dimension count 
we see that all the curves in $V$ with $\delta$ singularities have
precisely $\delta$ nodes.
\end{pf}

\begin{conj}\label{dconj}   
 $d_\delta(L)= T_\delta(L^2 , LK_S, K^2_S, c_2(S)).$
\end{conj}

Conjecture \ref{dconj} gives the hope of  
proving conjecture
\ref{mainconj}  via   the study of 
the cohomology of Hilbert schemes of points.

\begin{rem}
Note that one can generalize the above to singular points of arbitrary order:
Let $\mu=(m_1,\ldots,m_{l(\mu)})$  where $m_i\in\Z_{\ge 2}$. 
Let
$N(\mu):=\sum_{i=1}^s \binom{m_i+1}{2}$ and let 
  $S_\mu$ be the closure in $S^{[N(\mu)]}$ of the  
subset  of schemes of the form $\coprod_{i=1}^{l(\mu)}
\spec(\oo_{S,x_i}/\mmm_{S,x_i}^{m_i})$. 
Denote 
$
d_\mu(L):=\int_{S_\mu} c_{2l(\mu)}(L_{N(\mu)}).
$
We call a curve $D\in|L|$   of type $\mu$ if there are distinct points
$x_1,\ldots,x_{l(\mu)}$ in $S$, such that the ideal 
$\I_{D,x_i}$ is contained in $\mmm_{S,x_i}^{m_i}$. 
A straightforward generalization of the proof  of proposition \ref{port} shows
that, for $V$ a general $(N(\mu)-2l(\mu))$-dimensional linear subsystem of 
an $N(\mu)$-very ample line bundle $|L|$,
$d_\mu(L)$  counts the finite number of curves of type $\mu$ in $V$ 
with positive
multiplicities. Again I expect $d_\mu(L)$ to be a polynomial of degree 
$l(\mu)$
in $L^2$, $LK_S$, $K_S^2$ and $c_2(S)$.

In a similar way one can also deal with cusps instead of nodes.
\end{rem}


\begin{thebibliography}{ABCDEF}

\bibitem[B]{B} A. Beauville, {\em Counting curves on K3 surfaces}, preprint 
alg-geom/9701019.

\bibitem[Be]{Be} K. Behrend, 
{\em Gromov-Witten invariants in algebraic geometry},
Invent. Math. {\bf 127}, 601--617 (1997).

\bibitem[Be-F1]{Be-F1} K. Behrend, B. Fantechi, {\em
The intrinsic normal cone}, Invent Math. {\bf 128}, 45--88  (1997).



\bibitem[Be-F2]{Be-F2} K. Behrend, B. Fantechi, in preparation.


\bibitem[B-S]{B-S} M. Beltrametti, A.J. Sommese,  
{\em Zero cycles and $k$-th order
embeddings of smooth projective surfaces}, Problems on surfaces and their
classification, INDAM, Academic Press.

\bibitem[Bi-L]{Bi-L} Ch. Birkenhake, H. Lange,
{\em Complex Abelian varieties}, Grundlehren {\bf 302},
Springer Verlag Berlin Heidelberg 1992.

\bibitem[Br-Le]{Br-Le} J. Bryan, C. Leung, in preparation.

\bibitem[C]{C}  L. Caporaso, {\em Counting curves on surfaces: a guide
to new techniques and results}, preprint alg-geom/9611029.

\bibitem[C-H]{C-H} L. Caporaso, J. Harris, {\em Counting plane curves of 
any genus},  preprint alg-geom/9608025.

\bibitem[Ch]{Ch} Y. Choi, {\em Severi degrees in cogenus 4},
preprint alg-geom/9601013.

\bibitem[D-I]{D-I} P. Di Francesco, C. Itzykson,
{\em Quantum intersection rings}, The moduli space of curves, eds. 
R. Dijkgraaf,
C. Faber, G. van der Geer, Birkh\"auser, Boston 1995, pp 81--148.


\bibitem[F-G-vS]{F-G-vS} B. Fantechi, L. G\"ottsche, D. van Straten,
{\em  Euler number of the compactified Jacobian
and multiplicity of rational curves}, preprint 1997.

\bibitem[Fu]{Fu} W. Fulton, {\em Intersection theory}, Ergebnisse der 
Mathematik und ihrer Grenzgebiete, Springer-Verlag, Berlin  1984.

\bibitem[Ge]{Ge} E. Getzler, {\em Intersection theory on 
$\overline {\cal M}_{1,4}$ and elliptic Gromov-Witten invariants}, 
preprint alg-geom/9612004.



\bibitem[H]{H} J. Harris. {\em On the Severi variety}, Invent. Math. {\bf 84},
445--461 (1986).

\bibitem[H-P]{H-P} J. Harris, R. Pandharipande, 
{\em Severi degrees in Cogenus 3}, preprint 1995.

\bibitem[K-Z]{K-Z} Kaneko, D. Zagier, {\em A generalized Jacobi theta function
and quasimodular forms}, eds. R. Dijkgraaf, C. Faber, G. van der Geer,
Birkh\"auser, Boston 1995, pp 165--172.


\bibitem[K-P]{K-P} S. Kleiman, R. Piene, private communication.



\bibitem[K-M]{K-M} M. Kontsevich, Y. Manin, {\em 
Gromov-Witten classes, quantum cohomology and enumerative geometry},
Comm. Math. Phys. {\bf 164}, 525--562 (1994).

\bibitem[L-B]{L-B} H. Lange, C. Birkenhake,
{\em Complex Abelian Varieties}, Grundlehren der mathematischen Wissenschaften
302, Springer Verlag, Berlin 1992.



\bibitem[L-T]{L-T} J. Li, G. Tian,
{\em Virtual moduli cycles and GW-invariants}, preprint alg-geom/9602007.


\bibitem[R1]{R1} Z. Ran, {\em The degree of a Severi variety}, 
Bull. Amer. Math. Soc. {\bf 17}, (1987), 125--128.
\bibitem[R2]{R2} Z. Ran, {\em Enumerative geometry of singular 
plane curves}, Invent. Math. {\bf 97} 447--469 (1989).

\bibitem[Va]{Va} R. Vakil, {\em Counting curves of any genus
on rational ruled surfaces}, preprint 1997.


\bibitem[V1]{V1} I. Vainsencher, {\em Counting divisors with 
prescribed singularities}, Trans. Amer. Math. Soc.. {\bf 267}, 
399--422 (1981).


\bibitem[V2]{V2} I. Vainsencher, {\em Enumeration of $n$-fold tangent 
hyperplanes to a surface}, J. Alg. Geom.  , 503--526 (1995).
 
\bibitem[Y-Z]{Y-Z} S. T. Yau, E. Zaslow, {\em BPS states, string duality,
and nodal curves on K3}, preprint hep-th/9512121.
 
\end{thebibliography}
\end{document}